\documentclass[a4paper,12pt]{article}
\usepackage[dvips]{epsfig,graphicx}
\usepackage{amsmath,amsfonts,verbatim}

\newcommand{\RR}{{\mathbb R}}

\newcommand{\ths}{\theta^*}
\newcommand{\G}{\Gamma}

\newcommand{\h}{\hat}

\newcommand{\As}{A^*}
\newcommand{\Es}{E^*}
\newcommand{\mat}{{\rm Mat}}
\newcommand{\qed}{\hfill\hbox{\rule{3pt}{6pt}}}
\newcommand{\proof}{{\sc Proof. }}
\newcommand{\MX}{\mat_X(\RR)}

\newcommand{\CR}{{\cal R}}

\newtheorem{theorem}{Theorem}[section]
\newtheorem{lemma}[theorem]{Lemma}
\newtheorem{corollary}[theorem]{Corollary}
\newtheorem{proposition}[theorem]{Proposition}
\newtheorem{definition}[theorem]{Definition}
\newtheorem{note}[theorem]{Note}

\newtheorem{notation}[theorem]{Notation}

\title{BIPARTITE $Q$-POLYNOMIAL DISTANCE-REGULAR GRAPHS \\ AND UNIFORM POSETS}

\author{\v{S}tefko Miklavi\v{c} \\ 
        University of Primorska \\
        UP PINT and UP FAMNIT \\
        Muzejski trg 2 \\ 
        6000 Koper, Slovenia \\
        stefko.miklavic@upr.si \and
        Paul Terwilliger \\
        Department of Mathematics \\
        University of Wisconsin \\
        480 Lincoln Drive \\
        Madison WI 53706-1388, USA \\
        terwilli@math.wisc.edu}

\begin{document}
\maketitle
\begin{abstract}
  Let $\G$ denote a bipartite distance-regular graph with vertex set $X$ and diameter $D \ge 3$.
  Fix $x \in X$ and let $L$ (resp. $R$) denote the corresponding lowering (resp. raising) matrix.
  We show that each $Q$-polynomial structure for $\G$ yields a certain linear dependency among
  $RL^2$, $LRL$, $L^2R$, $L$. Define a partial order $\le$ on $X$ as follows. 
  For $y,z \in X$ let $y \le z$ whenever $\partial(x,y)+\partial(y,z)=\partial(x,z)$,
  where $\partial$ denotes path-length distance.
  We determine whether the above linear dependency gives this poset a 
  uniform  
   or strongly uniform structure.
  We show that except for one special 
  case a uniform structure is attained, and except for
 three special cases a  
  strongly uniform structure is attained.
\end{abstract}

\section{Introduction}
\label{sec:intro}

\noindent
In his thesis \cite{del} Delsarte introduced the $Q$-polynomial property for a 
distance-regular graph $\G$ (see Section \ref{sec:prelim} for formal definitions).
Since then the $Q$-polynomial property has been investigated by many authors,
such as Bannai and Ito \cite{BI}, Brouwer, Cohen and Neumaier \cite{BCN}, 
Caughman \cite{Caug1, Caug, Caug6, Caug5, Caug4, Caug3}, Curtin \cite{cur1, cur2}, Juri\v si\' c, Terwilliger, and \v Zitnik \cite{JTZ}, 
Lang \cite{lang1, lang2}, Lang and Terwilliger \cite{LT},
Miklavi\v c \cite{Mi1, Mi2, Mi3, Mi4}, Pascasio \cite{Pas1, Pas2}, Tanaka \cite{tanaka1, tanaka2},
Terwilliger \cite{ter4, ter1, ter6, ter5}, and Weng \cite{weng1, weng2}.

\medskip \noindent
To simplify this investigation, it is sometimes assumed that $\G$ is bipartite \cite{Caug1, Caug, Caug6, Caug5, Caug4, Caug3, lang1, lang2, Mi2, Mi3}
and this is the point of view taken in the present paper. For the rest of this Introduction assume
$\G$ is bipartite and $Q$-polynomial. To avoid trivialities assume $\G$ has diameter $D \ge 3$ and valency $k \ge 3$.

\medskip \noindent
In \cite{ter2} Terwilliger introduced the subconstituent algebra of $\G$. For each vertex $x$ of $\G$,
the corresponding subconstituent algebra $T=T(x)$ is generated by the adjacency matrix $A$ and a 
certain diagonal matrix $\As=\As(x)$. The eigenspaces of $\As$ are the subconstituents of $\G$ 
with respect to $x$. The matrices $A$ and $\As$ satisfy two relations called the 
tridiagonal relations \cite[Lemma 5.4]{ter3}, \cite{ter7}. 
The first (resp. second) tridiagonal relation is of degree $3$ in $A$ (resp. $\As$) and of degree $1$ in 
$\As$ (resp. $A$).
In \cite{ter3} the tridiagonal relations are used to describe the combinatorics of $\G$.
In this description it is natural to view $\G$ as the Hasse diagram for a 
ranked poset. The partial order $\le$ is defined as follows.
For vertices $y,z$ of $\G$ let $y \le z$ whenever $\partial(x,y)+\partial(y,z)=\partial(x,z)$,
where $\partial$ denotes path-length distance. 
The poset structure induces a decomposition
$A=L+R$, where $L=L(x)$ (resp. $R=R(x)$) is the lowering matrix (resp. raising matrix) of $\G$ with 
respect to $x$. For vertices $y,z$ of $\G$ the $(y,z)$-entry of $L$ is $1$ if $z$ covers $y$, and $0$ otherwise.
The matrix $R$ is the transpose of $L$.
In the first tridiagonal relation, if one eliminates $A$ using $A=L+R$, one finds that on each 
$x$-subconstituent of $\G$ the elements
$$
  RL^2, \quad LRL, \quad L^2R, \quad L
$$
are linearly dependent.
The coefficients in this linear dependence depend on the subconstituent. 
We call this collection of dependencies an $R/L$ dependency structure.

\medskip \noindent
Motivated by these $R/L$ dependency structures, in \cite{ter1} Terwilliger introduced the uniform property
for a partially ordered set. In that work he described the algebraic structure of the uniform posets
and displayed eleven infinite families of examples. 

\medskip \noindent
In spite of the known connection between the 
$Q$-polynomial property and uniform posets, a careful study of this connection was not completed
until now. The goal of the present paper is to provide this study.
As part of this study we introduce a variation on the 
uniform property called strongly uniform. Strongly uniform implies uniform.
For each $Q$-polynomial structure on $\G$ we determine precisely when the corresponding
$R/L$ dependency structure is uniform or
strongly uniform. 
To describe our results let $\{\theta_i\}_{i=0}^D$ denote 
the ordering of the eigenvalues of $\G$ for the given $Q$-polynomial structure.
Consider the following cases:
\begin{itemize}
\item[{\rm (i)}]
$\G$ is the hypercube $H(D,2)$ with $D$ even and
$\theta_i=(-1)^i (D-2i)$ for $0 \le i \le D$;
\item[{\rm (ii)}]
$\G$ is the antipodal quotient $\overline{H}(2D,2)$ and $\theta_i = 2D - 4i$ for $0 \le i \le D$;
\item[{\rm (iii)}]
$D=3$ and $\G$ is of McFarland type with parameters $(1,t)$ for some integer $t \ge 2$, and 
$\theta_0, \theta_1, \theta_2, \theta_3$ are $t(t+1), t, -t, -t(t+1)$ respectively.
\end{itemize}
(See Section \ref{sec:McFar} for the meaning of McFarland type).
In Case (i) the corresponding
$R/L$ dependency structure is not uniform.
In Cases (ii), (iii) this structure is uniform but not
strongly uniform.
In all other cases 
this structure is strongly uniform.

\medskip \noindent
The paper is organized as follows. In Sections \ref{sec:prelim} and \ref{sec:dualBM} we discuss the Bose-Mesner algebra
and the dual Bose-Mesner algebra of a distance-regular graph. In Sections \ref{sec:McFar} and \ref{sec:bip_poset} we consider the
bipartite case and discuss the associated poset structure. In Section \ref{sec:R/L} we consider $R/L$ dependency structures.
In Section \ref{sec:uniform} we review the uniform property and define the strongly uniform property. 
In Sections \ref{sec:beta2}--\ref{sec:qs1} we consider a given $Q$-polynomial structure for our graph. 
We determine precisely when the corresponding
$R/L$ dependency structure is uniform or strongly uniform.
Our main result is Theorem \ref{main:thm}.

\section{Preliminaries}
\label{sec:prelim}

Let $X$ denote a nonempty finite set. Let $\MX$ denote the $\RR$-algebra consisting of the matrices 
with entries in $\RR$, and rows and columns indexed by $X$. Let $V=\RR^X$ denote the vector space over 
$\RR$ consisting of the column vectors with entries in $\RR$ and rows indexed by $X$. Observe that $\MX$ 
acts on $V$ by left multiplication. We refer to $V$ as the {\em standard module} of $\MX$. 
We endow $V$ with the bilinear form $\langle \: , \rangle: V \times V \to \RR$ that satisfies $\langle u,v \rangle = u^t v$ for 
$u,v \in V$, where $t$ denotes transpose. For $y \in X$ let $\h{y}$ denote the vector in $V$ that has $y$-coordinate $1$ 
and all other coordinates $0$. Observe that $\{\hat{y}\,|\,y \in X\}$ is an orthonormal basis for $V$.

\medskip \noindent
Throughout the paper let $\G=(X,\CR)$ denote a finite, undirected, connected graph, without loops or multiple edges, 
with vertex set $X$, edge set $\CR$, path-length distance function $\partial$, and diameter 
$D:=\max \{\partial(x,y) \, | \, x,y \in X\}$. For $x \in X$ and an integer $i$ let
$\G_i(x) = \{ y \in X \, | \, \partial(x,y)=i \}$. We 
abbreviate $\G(x) = \G_1(x)$. For an integer $k \ge 0$ we say 
$\G$ is {\em regular with valency} $k$ whenever $|\G(x)|=k$ for all $x \in X$. We say $\G$ is 
{\em distance-regular} whenever for all integers $0 \le h,i,j \le D$ and all $x,y \in X$ with 
$\partial(x,y)=h$ the number $p_{ij}^h := | \G_i(x) \cap \G_j(y) |$
is independent of $x,y$. The constants $p_{ij}^h$ are known as the {\em intersection numbers} of $\G$.
For convenience set 
$c_i:=p_{1, i-1}^i \, (1 \le i \le D)$, $a_i:=p_{1i}^i \, (0 \le i \le D)$, 
$b_i:=p_{1, i+1}^i \, (0 \le i \le D-1)$, $k_i:=p_{ii}^0 \, (0 \le i \le D)$, and $c_0:=0$, $b_D:=0$.
For the rest of this paper assume $\Gamma$ is distance-regular with diameter $D \ge 3$.
By the triangle inequality, for $0 \le h,i,j \le D$ we have $p_{ij}^h = 0$ (resp. $p_{ij}^h \ne 0$) whenever one of $h,i,j$ is 
greater than (resp. equal to) the sum of the other two. 
In particular $c_i \ne 0$ for $1 \le i \le D$ and $b_i \ne 0$ for $0 \le i \le D-1$.
Observe that $\G$ is regular with valency $k=b_0=k_1$ and that 
$c_i+a_i+b_i=k$ for $0 \le i \le D$.

\medskip \noindent 
We recall the Bose-Mesner algebra of $\G$. For $0 \le i \le D$ let $A_i$ denote the matrix in 
$\MX$ with $(y,z)$-entry
\begin{equation}
\label{dm1}
  (A_i)_{y z} = \left\{ \begin{array}{ll}
                 1 & \hbox{if } \; \partial(y,z)=i, \\
                 0 & \hbox{if } \; \partial(y,z) \ne i \end{array} \right. \qquad (y,z \in X).
\end{equation}
We call $A_i$ the $i$th {\it distance matrix} of $\G$. We abbreviate $A:=A_1$ and call this the 
{\it adjacency matrix} of $\G.$ We observe
(ai)   $A_0 = I$;
(aii)  $J = \sum_{i=0}^D A_i$;
(aiii)  $A_i^t = A_i  \;(0 \le i \le D)$;
(aiv)   $A_iA_j = \sum_{h=0}^D p_{ij}^h A_h \;(0 \le i,j \le D)$,
where $I$ (resp. $J$) denotes the identity matrix (resp. all 1's matrix) in  $\MX$. Using these facts
we find $\{A_i\}_{i=0}^D$ is a basis for a commutative subalgebra $M$ of $\MX$. We call $M$ the 
{\it Bose-Mesner algebra} of $\G$. By \cite[p.~190]{BI} $A$ generates $M$. By 
\cite[p.~45]{BCN} $M$ has a basis $\{E_i\}_{i=0}^D$ such that 
(ei)   $E_0 = |X|^{-1}J$;
(eii)  $I = \sum_{i=0}^D E_i$;
(eiii)  $E_i^t =E_i  \;(0 \le i \le D)$;
(eiv)   $E_iE_j =\delta_{ij}E_i  \;(0 \le i,j \le D)$.
We call $\{E_i\}_{i=0}^D$  the {\it primitive idempotents} of $\Gamma$.  
The primitive idempotent $E_0$ is said to be {\em trivial}. 

\noindent
We recall the eigenvalues of $\G$. Since $\{E_i\}_{i=0}^D$ form a basis for $M$,
there exist scalars $\{\theta_i\}_{i=0}^D$ in $\RR$ such that $A = \sum_{i=0}^D \theta_i E_i.$
Combining this with (eiv) we find
$$
  A E_i = E_i A = \theta_i E_i \qquad \qquad (0 \le i \le D).
$$ 
We call $\theta_i$ the {\em eigenvalue} of $\G$ associated with $E_i$. 
The $\{\theta_i\}_{i=0}^D$ are mutually distinct since $A$ generates $M$. 
By (ei) we have $\theta_0 = k$. 
By (eii)--(eiv), 
\begin{equation}
\label{de}
  V = E_0V + E_1V + \cdots + E_DV \qquad \hbox{(orthogonal direct sum)}.
\end{equation}
For $0 \le i \le D$ the space $E_iV$ is the eigenspace of $A$ associated with $\theta_i$.
Let $m_i$ denote the rank of $E_i$ and note that $m_i$ is the dimension of $E_iV$.
We call $m_i$ the {\em multiplicity} of $\theta_i$. 

\medskip \noindent
We recall the Krein parameters of $\G$.
Let $\circ $ denote the entrywise product in $\MX$. Observe that
$A_i\circ A_j= \delta_{ij}A_i$ for $0 \leq i,j\leq D$, so $M$ is closed under $\circ$. Thus there 
exist scalars $q^h_{ij} \in \RR$  $(0 \leq h,i,j\leq D)$ such that
$$
  E_i \circ E_j = |X|^{-1}\sum_{h=0}^D q^h_{ij}E_h \qquad \qquad (0 \leq i,j\leq D).
$$
The parameters $q^h_{ij}$ are called the {\it Krein parameters of} $\G$.
By \cite[Proposition 4.1.5]{BCN} these parameters are nonnegative.
The given ordering $\{E_i\}_{i=0}^D$ of the primitive idempotents
is said to be $Q$-{\em polynomial} whenever for $0 \le h,i,j \le D$ 
the Krein parameter $q^h_{ij}=0$ (resp. $q^h_{ij} \ne 0$) whenever one of $h,i,j$ is greater than (resp. equal to) 
the sum of the other two. Let $E$ denote a nontrivial primitive idempotent of $\G$ and let $\theta$ denote the 
corresponding eigenvalue. We say that $\G$ is {\em Q-polynomial with respect to E} (or $\theta$) whenever there exists a 
$Q$-polynomial ordering $\{E_i\}_{i=0}^D$ of the primitive idempotents of $\G$ such that $E_1 = E$. 

\section{The dual Bose-Mesner algebra}
\label{sec:dualBM}

We continue to discuss the distance-regular graph $\G$ from Section \ref{sec:prelim}.
In this section we recall the dual Bose-Mesner algebra of $\G$. For the rest of the paper fix $x \in X$. 
For $ 0 \le i \le D$ let $E_i^*=E_i^*(x)$ denote the diagonal matrix in $\MX$ with $(y,y)$-entry 
\begin{equation}
\label{den0}
(\Es_i)_{y y} = \left\{ \begin{array}{lll}
                 1 & \hbox{if } \; \partial(x,y)=i, \\
                 0 & \hbox{if } \; \partial(x,y) \ne i \end{array} \right. 
                 \qquad (y \in X).
\end{equation}
We call $E_i^*$ the  $i$th {\it dual idempotent} of $\G$ with respect to $x$ \cite[p.~378]{ter2}. 
For convenience set $\Es_i = 0$ for $i < 0$ or $i > D$.
We observe
(esi) $I = \sum_{i=0}^D \Es_i$;
(esii) $E_i^{*t} = E_i^* \;\; (0 \le i \le D)$;
(esiii) $E_i^* E_j^* =  \delta_{ij} E_i^*\;\; (0 \le i,j \le D)$.
By these facts $\{\Es_i\}_{i=0}^D$ form a basis for a commutative subalgebra $M^*=M^*(x)$ 
of $\MX$. We call $M^*$ the {\it dual Bose-Mesner algebra} of $\G$ with respect to 
$x$ \cite[p.~378]{ter2}. By (esi)--(esiii),
\begin{equation}
\label{vsub}
V = E_0^*V + E_1^*V + \cdots + E_D^*V \qquad \qquad {\rm (orthogonal\ direct\ sum}).
\end{equation}
For $0 \le i \le D$ the subspace $\Es_iV$ has basis $\{\hat{y} \, | \, y \in \G_i(x)\}$.
Moreover the dimension of $\Es_iV$ is $k_i$.

\medskip \noindent
The algebras $M$ and $M^*$ are related as follows.
By \cite[Lemma 3.2]{ter2},
\begin{equation}
\label{alg_rel}
  \Es_i A_j \Es_h = 0 \;\; \hbox{if and only if} \; \; p_{ij}^h=0 \qquad (0 \le h, i,j \le D).
\end{equation}
Let $E$ denote a nontrivial primitive idempotent of $\G$ and assume $\G$ is $Q$-polynomial with respect to $E$.
Let $A^*=A^*(x)$ denote the diagonal matrix in $\MX$ with $(y,y)$-entry
$$
  A^*_{yy} = |X| E_{xy} \qquad (y \in X).
$$
We call $A^*$ the {\em dual adjacency matrix of} $\G$ that corresponds to $E$ and $x$. 
By \cite[Lemma 3.11(ii)]{ter2} $A^*$ generates $M^*$. We recall the dual eigenvalues for our $Q$-polynomial structure.
Since $\{\Es_i\}_{i=0}^D$ form a basis for $M^*$ there exist scalars $\{\ths_i\}_{i=0}^D$ in $\RR$ 
such that $\As = \sum_{i=0}^D \ths_i E_i^*$. Combining this with (esiii) we find 
\begin{equation}
\label{qpoly:eq1}
\As E_i^* = E_i^* \As =\ths_i E_i^* \qquad \qquad (0 \le i \le D).
\end{equation} 
We call $\{\ths_i\}_{i=0}^D$ the {\em dual eigenvalue sequence} 
for the given $Q$-polynomial structure.
The $\{\ths_i\}_{i=0}^D$ are mutually distinct since $\As$ generates $M^*$. 
For $0 \le i \le D$ the space $\Es_iV$ is the eigenspace of $\As$ associated with $\ths_i$.
By \cite[Proposition 3.4.(iv)]{BI} we have that $\ths_0 = {\rm rank}(E)$.
Let $\theta$ denote the eigenvalue of $\G$ associated with $E$. By \cite[p. 128]{BCN},
\begin{equation}
\label{eq:de}
  c_i \ths_{i-1} + a_i \ths_i +  b_i \ths_{i+1} = \theta \ths_i \qquad (0 \le i \le D),
\end{equation}
where $\ths_{-1}$ and $\ths_{D+1}$ are indeterminants.

\begin{lemma}
\label{lem:fundamental1}
{\rm (\cite[Lemma 5.4]{ter3})} 
Let $\{E_i\}_{i=0}^D$ denote a $Q$-polynomial ordering of the primitive idempotents 
of $\G$ and for $0 \le i \le D$ let $\theta_i$ denote the eigenvalue of $\G$ for $E_i$. 
Let $\{\ths_i\}_{i=0}^D$ denote the dual eigenvalue sequence for the given 
$Q$-polynomial structure. Then the following {\rm (i)--(iii)} hold.
\begin{itemize}
\item[{\rm (i)}]
There exists $\beta \in \RR$ such that 
\begin{equation}
\label{eq:beta}
  \beta +1 = \frac{\theta_{i-2} - \theta_{i+1}}{\theta_{i-1} - \theta_i}
           = \frac{\ths_{i-2} - \ths_{i+1}}{\ths_{i-1} - \ths_i}
\end{equation}
for $2 \le i \le D-1$.
\item[{\rm (ii)}]
There exist $\gamma, \gamma^* \in \RR$ such that both
\begin{equation}
\label{eq:gamma}
  \gamma = \theta_{i-1} - \beta \theta_i + \theta_{i+1}, \qquad \qquad \gamma^* = \ths_{i-1} - \beta \ths_i + \ths_{i+1}
\end{equation}
for $1 \le i \le D-1$.
\item[{\rm (iii)}]
There exist $\varrho, \varrho^* \in \RR$ such that both
\begin{equation}
\label{eq:rho}
\begin{split}
  \varrho & = \theta_{i-1}^2 - \beta \theta_{i-1} \theta_i + \theta_i^2 - \gamma (\theta_{i-1} + \theta_i), \\ 
  \varrho^* & = \ths_{i-1} \hspace{-4mm}\phantom{.}^2 - \beta \ths_{i-1} \ths_i + \ths_i \hspace{-1mm}\phantom{.}^2 - \gamma^* (\ths_{i-1} + \ths_i)
\end{split}
\end{equation}
for $1 \le i \le D$.
\end{itemize}
\end{lemma} 

\begin{lemma}
\label{lem:fundamental2}
{\rm (\cite[Lemma 5.4]{ter3})}
Let $E$ denote a $Q$-polynomial primitive idempotent of $\G$ and 
let $\As = \As(x)$ denote the corresponding dual adjacency matrix.
Then both 
\begin{equation}
\label{fund_eq}
  [A, A^2 \As -\beta A \As A + \As A^2 - \gamma ( A \As + \As A)- \varrho \As] = 0, 
\end{equation}
\begin{equation}
  [\As, \As \hspace{-1mm}\phantom{.}^2 A -\beta \As A \As + A \As \hspace{-1mm} \phantom{.}^2 - \gamma^* ( \As A + A \As) - \varrho^* A] = 0,
\end{equation}
where $[r,s] = rs-sr$ and $\beta, \gamma, \gamma^*, \varrho, \varrho^*$ are from Lemma \ref{lem:fundamental1}.
\end{lemma}

\section{Bipartite distance-regular graphs}
\label{sec:McFar}

We continue to discuss the distance-regular graph $\G$ from Section \ref{sec:prelim}.
Recall that $\G$ is {\em bipartite} whenever $a_i=0$ for $0 \le i \le D$. 
For $\G$ bipartite, $p^h_{ij}=0$ if $h+i+j$ is odd $(0 \le h,i,j \le D)$.
In this case 
\begin{equation}
\label{prelim:eq1}
  \Es_i A \Es_h = 0 \quad \hbox{if} \quad |h-i| \ne 1 \qquad (0 \le h,i \le D).
\end{equation}
%

\medskip \noindent
The case in which $\G$ is bipartite with $D=3$ will play an important role.

\medskip \noindent
By \cite[Theorem 1.6.1.]{BCN}, $\G$ is bipartite with $D=3$ 
if and only if $\G$ is the incidence graph of a square $2$-$(v,k,\lambda)$ design. 
In this case $c_2 = \lambda$ and $v=1+k(k-1)/\lambda$. 
See \cite{BJL} for more information and background on square $2$-designs.

\medskip \noindent
Pick integers $d \ge 1$ and $t \ge 2$. A square $2$-$(v,k,\lambda)$ design is said to have
{\em McFarland type with parameters} $(d,t)$ whenever 
$$
  v = t^{d+1} \Bigg( 1 + \frac{t^{d+1}-1}{t-1} \Bigg), \quad k = t^d \: \frac{t^{d+1}-1}{t-1}, \quad \lambda = t^d \: \frac{t^d-1}{t-1}.
$$
For the moment assume that $t$ is a prime power. By \cite[Corollary II.8.17]{BJL}, a square $2$-design of McFarland type with parameters $(d,t)$
exists for every integer $d \ge 1$. By \cite[p. 982]{BJL} this design can be realized as a McFarland difference set. 

\medskip \noindent
Our graph $\G$ is said to have {\em McFarland type with parameters} $(d,t)$ whenever $\G$ is 
the incidence graph of a square $2$-design of McFarland type with parameters $(d,t)$.

\section{The bipartite case; lowering and raising matrices}
\label{sec:bip_poset}

We continue to discuss the distance-regular graph $\G$ from Section \ref{sec:prelim}.
For the rest of this paper assume that $\G$ is bipartite. 

\medskip \noindent
Define a partial order $\le$ on $X$ such that for all $y,z \in X$,
$$
  y \le z \quad \hbox{if and only if} \quad \partial(x,y)+\partial(y,z)=\partial(x,z).
$$
For $y,z \in X$ define $y < z$ whenever $y \le z$ and $y \ne z$. 
We say that $z$ {\em covers} $y$ whenever $y < z$ and there does not exist $w \in X$ such that $y < w < z$.
Note that $z$ covers $y$ if and only if $y,z$ are adjacent and $\partial(x,y)+1=\partial(x,z)$.
For $0 \le i \le D$ each vertex in $\G_i(x)$ covers exactly $c_i$ vertices from $\G_{i-1}(x)$, and is covered by 
exactly $b_i$ vertices in $\G_{i+1}(x)$. Therefore the partition $\{\G_i(x)\}_{i=0}^D$ of $X$ 
is a {\em grading} of the poset $(X,\le)$ in the sense of \cite[Section 1]{ter1}. 

\begin{definition}
\label{def:LR}
{\rm 
Define matrices $L=L(x)$ and $R=R(x)$ by
$$
  L = \sum_{i=1}^D \Es_{i-1} A \Es_i, \qquad \qquad
  R = \sum_{i=0}^{D-1} \Es_{i+1} A \Es_i.
$$
Note that $R=L^t$ and $L+R=A$.}
\end{definition}
We have three observations.

\begin{lemma}
\label{lem:rl-covers}
Let $L,R$ be as in Definition \ref{def:LR}.
Then the following {\rm (i), (ii)} hold for $y \in X$.
\begin{itemize}
\item[{\rm (i)}]
$L \h{y} = \sum \h{z}$, where the sum is over all $z \in X$ that are covered by $y$. 
\item[{\rm (ii)}]
$R \h{y} = \sum \h{z}$, where the sum is over all $z \in X$ that cover $y$.
\end{itemize}
\end{lemma}
Motivated by Lemma \ref{lem:rl-covers} we call $L$ (resp. $R$) the {\em lowering matrix}
(resp. {\em raising matrix}) {\em of} $\G$ {\em with respect to} $x$.

\begin{lemma}
\label{lem:rl-covers2}
Let $L,R$ be as in Definition \ref{def:LR}.
Then the following {\rm (i), (ii)} hold.
\begin{itemize}
\item[{\rm (i)}]
$R \Es_i V \subseteq \Es_{i+1}V$ for $0 \le i \le D-1$, and $R \Es_D V = 0$;
\item[{\rm (ii)}]
$L \Es_i V \subseteq \Es_{i-1}V$ for $1 \le i \le D$, and $L \Es_0 V = 0$.
\end{itemize}
\end{lemma}

\begin{lemma}
\label{main:lem1}
Let $L,R$ be as in Definition \ref{def:LR}. Then for $1 \le i \le D$ the following {\rm (i)--(iv)} hold.
\begin{itemize}
\item[{\rm (i)}]
$\Es_{i-1} A \Es_i = L \Es_i$;
\item[{\rm (ii)}]
$\Es_{i-1} A \Es_i = \Es_{i-1} L$;
\item[{\rm (iii)}]
$\Es_i A \Es_{i-1} = R \Es_{i-1}$;
\item[{\rm (iv)}]
$\Es_i A \Es_{i-1} = \Es_i R$.
\end{itemize}
Moreover
\begin{equation}
\label{eq_pom1}
  L \Es_0 = 0, \qquad \Es_D L = 0, \qquad R \Es_D = 0, \qquad \Es_0 R = 0.
\end{equation}
\end{lemma} 

\begin{lemma}
\label{lem:rl}
Let $L,R$ be as in Definition \ref{def:LR}.
Then the following {\rm (i)--(iii)} hold for $0 \le i \le D$.
\begin{itemize}
\item[{\rm (i)}]
$R L^2 \Es_i = \Es_{i-1} A \Es_{i-2} A \Es_{i-1} A \Es_i$;
\item[{\rm (ii)}]
$L R L \Es_i = \Es_{i-1} A \Es_i A \Es_{i-1} A \Es_i$;
\item[{\rm (iii)}]
$L^2 R \Es_i = \Es_{i-1} A \Es_i A \Es_{i+1} A \Es_i$.
\end{itemize}
\end{lemma}
\proof
Immediate from (esiii) and Definition \ref{def:LR}. \qed 

%

\begin{lemma}
\label{main:lem4}
Let $L,R$ be as in Definition \ref{def:LR}.
Then 
$$
  \Es_{i-1} A^3 \Es_i = R L^2 \Es_i + L R L \Es_i + L^2 R \Es_i
$$ 
for $1 \le i \le D$.
\end{lemma}
\proof
Straightforward using $A=L+R$ and Lemma \ref{lem:rl-covers2}. \qed

\medskip \noindent
From now on we use the following notational convention.
\begin{notation}
\label{blank1}
{\rm For the rest of this paper we assume our distance-regular graph $\G$ is bipartite with valency $k \ge 3$.
Let $\{E_i\}_{i=0}^D$ denote a $Q$-polynomial ordering of the primitive idempotents of $\G$ and let
$\{\theta_i\}_{i=0}^D$ denote the corresponding eigenvalues.
Abbreviate $E=E_1$.
Recall our fixed vertex $x \in X$ from Section \ref{sec:dualBM}. 
For $0 \le i \le D$ let
$\Es_i = \Es_i(x)$ denote the $i$th dual idempotent of $\G$ with respect to $x$.
Let $\As = \As(x)$ denote the dual adjacency matrix of $\G$ that corresponds to $E$ and $x$. Let
$\{\ths_i\}_{i=0}^D$ denote the dual eigenvalue sequence for the given $Q$-polynomial structure.
Let the scalars $\beta, \gamma, \gamma^*, \varrho, \varrho^*$ be from Lemma \ref{lem:fundamental1}. 
Let the matrices $L=L(x)$ and $R=R(x)$ be as in Definition \ref{def:LR}.}
\end{notation}

\noindent
With reference to Notation \ref{blank1}, we have $\gamma = 0$ by \cite[Theorem 8.2.1]{BCN} and since $\G$ is bipartite. 
Thus by \eqref{fund_eq},
\begin{equation}
\label{fund_eq_bip}
  [A, A^2 \As -\beta A \As A + \As A^2 - \varrho \As] = 0.
\end{equation}

\section{The $R/L$ dependency structure}
\label{sec:R/L}

In this section we display certain linear dependencies among $RL^2, RLR, L^2R,L$.

\begin{lemma}
\label{main:lem3}
With reference to Notation \ref{blank1} the following {\rm (i), (ii)} hold for $1 \le i \le D$.
\begin{itemize} 
\item[{\rm (i)}]
$\Es_{i-1} A^2 \As A \Es_i = \ths_{i-1} R L^2 \Es_i + \ths_{i-1} L R L \Es_i + \ths_{i+1} L^2 R \Es_i$;
\item[{\rm (ii)}]
$\Es_{i-1} A \As A^2 \Es_i = \ths_{i-2} R L^2 \Es_i + \ths_i L R L \Es_i + \ths_i L^2 R \Es_i$.
\end{itemize}
\end{lemma} 
\proof
Straightforward using $A=L+R$ along with \eqref{qpoly:eq1} and Lemma \ref{lem:rl-covers2}. \qed

\begin{proposition}
\label{main:prop5a}
With reference to Notation \ref{blank1}, for $1 \le i \le D$ the equation 
\begin{equation}
\label{general}
\begin{split}
  \frac{\ths_i-\ths_{i-1} + (\beta+1) (\ths_{i-2}-\ths_{i-1})}{\ths_i - \ths_{i-1}} & R L^2 + 
  (\beta+2) L R L \\ 
  &+ \frac{\ths_i - \ths_{i-1} + (\beta+1) (\ths_i - \ths_{i+1})}{{\ths_i - \ths_{i-1}}} L^2 R = \varrho  L
\end{split}
\end{equation}
holds on $\Es_iV$.
\end{proposition}
\proof
Multiply \eqref{fund_eq_bip} by $\Es_{i-1}$ on the left and by $\Es_i$ on the right. 
Divide the result by $\ths_{i-1}-\ths_i$ and simplify
using \eqref{qpoly:eq1} along with 
Lemmas \ref{main:lem1}(i), \ref{main:lem4}, \ref{main:lem3}. \qed

\medskip \noindent
We call the equations \eqref{general} the {\em R/L dependency structure} that corresponds to the given $Q$-polynomial structure.
We have a comment about the coefficients in line \eqref{general}.

\begin{lemma}
\label{main:lem6}
With reference to Notation \ref{blank1} the following {\rm (i), (ii)} hold.
\begin{itemize}
\item[{\rm (i)}]
For $3 \le i \le D$,
$$
  \frac{\ths_i - \ths_{i-1} + (\beta+1) (\ths_{i-2} - \ths_{i-1})}{\ths_i - \ths_{i-1}} = \frac{\ths_{i-3} - \ths_{i-1}}{\ths_i - \ths_{i-1}}.
$$
\item[{\rm (ii)}]
For $1 \le i \le D-2$,
$$
  \frac{\ths_i - \ths_{i-1} + (\beta+1) (\ths_i - \ths_{i+1})}{\ths_i - \ths_{i-1}} = \frac{\ths_i-\ths_{i+2}}{\ths_i - \ths_{i-1}}.
$$
\end{itemize}
\end{lemma}
\proof
(i) Evaluate the left-hand side using $\beta + 1 = (\ths_{i-3} - \ths_i)/(\ths_{i-2} - \ths_{i-1})$.

\smallskip \noindent
(ii) Evaluate the left-hand side using $\beta + 1 = (\ths_{i-1} - \ths_{i+2})/(\ths_i - \ths_{i+1})$.
\qed

\section{Uniform structures on a poset}
\label{sec:uniform}

In this section we discuss the uniform property for a partially ordered set \cite{ter1}. 
This property involves the notion of a parameter matrix.
With reference to Notation \ref{blank1}, by a
{\em parameter matrix} we mean a tridiagonal matrix 
$U=(e_{ij})_{1 \le i,j \le D}$ with entries in $\RR$ such that 
\begin{itemize}
\item[(1)]
$e_{ii}=1$ for $1 \le i \le D$;
\item[(2)]
$e_{i,i-1} \ne 0$  for
$2 \le i \le D$ or
 $e_{i-1,i} \ne 0$ for $2 \le i \le D$;
\item[(3)]
the principal submatrix 
$(e_{ij})_{r \le i,j \le p}$ is nonsingular
for $1 \le r \le p \le D$.
\end{itemize}
We abbreviate $e_i^- := e_{i,i-1}$ for $2 \le i \le D$ and 
$e_i^+ := e_{i,i+1}$ for $1 \le i \le D-1$.
For notational convenience define $e_1^- := 0$ and $e_D^+ := 0$.

\medskip \noindent 
By a {\em uniform structure} for $\G$ we mean a pair $(U,f)$ where  
$U=(e_{ij})_{1 \le i,j \le D}$ is a parameter matrix and $f=\{f_i\}_{i=1}^D$ is a vector in $\RR^D$ such that the equation
\begin{equation}
\label{eq:uniform}
  e_i^- R L^2 + L R L + e_i^+ L^2 R = f_i L
\end{equation}
holds on $\Es_iV$ for $1 \le i \le D$. 
By a {\em strongly uniform structure} for $\G$
we mean a uniform structure $(U,f)$ for $\G$ such that
$e_{i,i-1} \ne 0$ and $e_{i-1,i} \ne 0$ for $2 \le i \le D$. 
Note that a strongly uniform structure is uniform.


\begin{lemma}
\label{lem:transp}
With reference to Notation \ref{blank1} let $(U,f)$ denote a uniform structure for $\G$. 
Then the equation
$$
  e_i^- R^2 L + R L R + e_i^+ L R^2 = f_i R
$$
holds on $\Es_{i-1}V$ for $1 \le i \le D$.
\end{lemma}
\proof
The equation
(\ref{eq:uniform})
holds on $\Es_iV$ so
\begin{equation}
\label{eq_pom}
  (e_i^- R L^2 + L R L + e_i^+ L^2 R - f_i L)\Es_i=0.
\end{equation}
By Lemma \ref{main:lem1} we have $L \Es_j = \Es_{j-1} L$ and 
$\Es_jR = R\Es_{j-1}$ for $1 \le j \le D$. Evaluating \eqref{eq_pom} using this and \eqref{eq_pom1}
we find 
\begin{equation}
\label{eq_pom2}
  \Es_{i-1}(e_i^- R L^2 + L R L + e_i^+ L^2 R - f_i L)=0.
\end{equation}
In line \eqref{eq_pom2} apply the transpose map to each term and
recall $R=L^t$. This yields
$$
(e_i^- R^2 L + R L R + e_i^+ L R^2 - f_i R)\Es_{i-1}=0
$$
and
the result follows. \qed

\medskip \noindent
See \cite{ter1} for more information on uniform posets.

\medskip \noindent
Recall our $Q$-polynomial structure from Notation \ref{blank1}.
Our next goal is to determine in which cases the corresponding $R/L$ dependency structure
is uniform or strongly uniform.
We first consider the case in which $\beta=-2$, where $\beta$ is from line \eqref{eq:beta}.

\section{The case $\mathbf{\beta=-2}$}
\label{sec:beta2}

Recall our $Q$-polynomial structure from Notation \ref{blank1}. In this section we determine whether
the corresponding $R/L$ dependency structure is  uniform or strongly uniform,
for the case $\beta = -2$.
We will be discussing the $D$-dimensional hypercube $H(D,2)$.
By \cite[Theorem 9.2.1]{BCN} $H(D,2)$ is distance-regular with diameter $D$ and intersection numbers
\begin{equation}
\label{int_num}
  b_i = D-i, \quad c_i=i \qquad (0 \le i \le D).
\end{equation}
By \cite[Theorem 9.2.1]{BCN} the eigenvalues of $H(D,2)$ are $\{D - 2i\}_{i=0}^D$.
By \cite[p.~304]{BI} the ordering $\{D - 2i\}_{i=0}^D$ is $Q$-polynomial.
For this $Q$-polynomial structure $\beta=2$.
If $D$ is odd then this $Q$-polynomial structure is unique. 
If $D$ is even then $H(D,2)$ has exactly one more $Q$-polynomial structure,
with eigenvalue ordering
$\{(-1)^i (D-2i)\}_{i=0}^D$ \cite[p. 305]{BI}. For this $Q$-polynomial structure $\beta=-2$.

\begin{proposition}
\label{prop:beta-2}
{\rm (\cite[Theorem 2]{ter4})}
With reference to Notation \ref{blank1} assume $\beta=-2$.
Then $D$ is even and $\G$ is $H(D,2)$ with 
the following $Q$-polynomial ordering of the eigenvalues: 
\begin{equation}
\label{eig_ord}
  \theta_i = (-1)^i (D-2i) \qquad (0 \le i \le D).
\end{equation} 
\end{proposition}

\begin{lemma}
\label{lem:beta-2}
With reference to Notation \ref{blank1}, assume $\G$ is $H(D,2)$ and let 
$\{\theta_i\}_{i=0}^D$ denote the $Q$-polynomial ordering of the eigenvalues \eqref{eig_ord}.
Let $\{\ths_i\}_{i=0}^D$ denote the corresponding dual eigenvalue sequence.
Then $\ths_i = \theta_i$ for $0 \le i \le D$.
Also 
$$
  \beta=-2, \qquad \gamma^*=0, \qquad \varrho=4, \qquad \varrho^*=4.
$$
\end{lemma}
\proof
We have $\ths_i = \theta_i$ by \cite[Theorem 1.1]{Caug1}. 
The remaining assertions follow from Lemma \ref{lem:fundamental1}. \qed

\begin{proposition}
\label{prop:beta-2a}
With reference to Notation \ref{blank1}, assume $\G$ is $H(D,2)$ and consider the $Q$-polynomial ordering of 
the eigenvalues \eqref{eig_ord}.
Then the corresponding $R/L$ dependency structure is that the equation
\begin{equation}
\label{eq:cube}
 \frac{2}{2i-D-1}R L^2 - \frac{2}{2i-D-1} L^2 R = 4 L
\end{equation}
holds on $\Es_iV$ for $1 \le i \le D$.
\end{proposition}
\proof 
Evaluate \eqref{general} using Lemma \ref{lem:beta-2}. \qed

\begin{proposition}
\label{pro:beta-2}
With reference to Notation \ref{blank1}, assume $\G$ is $H(D,2)$ and consider the $Q$-polynomial ordering of 
the eigenvalues \eqref{eig_ord}.
Then the corresponding $R/L$ dependency structure is not uniform.
\end{proposition}
\proof
The equation \eqref{eq:cube} does not match the form \eqref{eq:uniform}. \qed

\section{The case $\mathbf{\beta \ne -2}$}
\label{sec:beta_ne_-2}

Recall our $Q$-polynomial structure from Notation \ref{blank1}. 
Until further notice assume $\beta \ne -2$.
Under this assumption we show that the corresponding 
$R/L$ dependency structure is uniform.
Moreover we show that this structure is strongly uniform except in two special cases.
The following definition is for notational convenience.

\begin{definition}
\label{def:ef}
{\rm With reference to Notation \ref{blank1} assume $\beta \ne -2$. 
Let $U=(e_{ij})_{1 \le i,j \le D}$ denote the tridiagonal matrix with entries
\begin{align*}
e_{ii}   &= 1 &(1 \le i \le D), \\
e_{i,i-1} &= \frac{\ths_i - \ths_{i-1} + (\beta+1) (\ths_{i-2} - \ths_{i-1})}{(\beta+2)(\ths_i - \ths_{i-1})} \qquad &(2 \le i \le D), \\
e_{i-1,i} &= \frac{\ths_{i-1} - \ths_{i-2} + (\beta+1) (\ths_{i-1} - \ths_i)}{(\beta+2)(\ths_{i-1} - \ths_{i-2})} &(2 \le i \le D).
\end{align*}
For notational convenience write $e_i^- = e_{i,i-1}$ for $2 \le i \le D$ and $e_1^- = 0$, and also
$e_i^+ = e_{i,i+1}$ for $1 \le i \le D-1$ and $e_D^+ = 0$.
Define a vector $\{f_i\}_{i=1}^D$ in $\RR^D$ such that $f_i = \varrho/(\beta+2)$ for $1 \le i \le D$.}
\end{definition}

\begin{corollary}
\label{main:cor6} 
With reference to Notation \ref{blank1} and Definition \ref{def:ef} the following {\rm (i), (ii)} hold.
\begin{itemize}
\item[{\rm (i)}] For $3 \le i \le D$,
$$
  e_i^- = \frac{\ths_{i-3} - \ths_{i-1}}{(\beta+2)(\ths_i - \ths_{i-1})}.
$$
\item[{\rm (ii)}] For $1 \le i \le D-2$, 
$$
  e_i^+ = \frac{\ths_i-\ths_{i+2}}{(\beta+2)(\ths_i - \ths_{i-1})}.
$$
\end{itemize}
\end{corollary}
\proof
Immediate from Lemma \ref{main:lem6} and Definition \ref{def:ef}. \qed

\begin{proposition}
\label{main:prop5}
With reference to Notation \ref{blank1} and Definition \ref{def:ef}, the equation
\begin{equation}
\label{eq:uniform1}
  e_i^- RL^2 + LRL + e_i^+ L^2R = f_i L 
\end{equation}
holds on $\Es_iV$ for $1 \le i \le D$.
\end{proposition}
\proof
Divide \eqref{general} by $\beta+2$ and use Definition \ref{def:ef}. \qed

\medskip \noindent
Our next general goal is to determine whether the equations \eqref{eq:uniform1}
give a uniform or strongly uniform structure.
In order to do this we introduce some parameters $q$ and $s^*$. 

\section{The parameters $q$ and $s^*$}
\label{sec:qs}

Recall our $Q$-polynomial structure from Notation \ref{blank1}.
We would like to write the corresponding data in terms of two parameters 
$q$ and $s^*$. However, it will be convenient to exclude several special cases. 
The first special case is $H(D,2)$ with eigenvalue ordering $\{D-2i\}_{i=0}^D$.
The next special case concerns the antipodal quotient of $H(2D,2)$.
We denote this quotient graph by $\overline{H}(2D,2)$.
By \cite[p. 264]{BCN} $\overline{H}(2D,2)$ is distance-regular with diameter $D$ and intersection numbers 
$$
  b_i = 2D-i, \quad c_i = i \qquad (0 \le i \le D-1)
$$
and $c_D=2D$. By \cite[p. 264]{BCN} the eigenvalues of $\overline{H}(2D,2)$ are
\begin{equation}
\label{eig_aqhyp}
  \theta_i = 2D-4i \qquad (0 \le i \le D).
\end{equation}
By \cite[p. 306]{BI} the ordering \eqref{eig_aqhyp} is the unique $Q$-polynomial structure for $\overline{H}(2D,2)$. 
In order to describe some more special cases, we turn our attention to Notation \ref{blank1} with $D=3$.
By \cite[Proposition 4.2.2.(ii)]{BCN}, $b_2=1$ if and only if $\G$ is antipodal. 
In this case $b_1=k-1$, $c_2=k-1$, $c_3=k$. Moreover $\G$ has a unique 
$Q$-polynomial structure with eigenvalues $\theta_0=k$, $\theta_1=1$, $\theta_2=-1$, $\theta_3=-k$ \cite[p. 432]{BCN}.
For $b_2 > 1$, $\G$ has exactly two $Q$-polynomial structures:
$\theta_0=k$, $\theta_1=\sqrt{b_2}$, $\theta_2=-\sqrt{b_2}$, $\theta_3=-k$ and
$\theta_0=k$, $\theta_1=-\sqrt{b_2}$, $\theta_2=\sqrt{b_2}$, $\theta_3=-k$ \cite[p. 432]{BCN}.
In the following table we summarize the cases discussed so far.

\renewcommand\arraystretch{2}

\begin{table}[!h]
\begin{center}
\begin{tabular}{|c|c|c|}
 \hline 
 Case & $\G$ & $Q$-polynomial structure \\ \hline \hline
 I    & $H(D,2)$ & $\theta_i = D-2i \quad (0 \le i \le D)$ \\ \hline
 II   & $\overline{H}(2D,2)$ & $\theta_i = 2D - 4i \quad (0 \le i \le D)$ \\ \hline
 III  & $D=3$, $b_2=1$ & $\theta_0 = k$, $\theta_1 = 1$, $\theta_2 = -1$, $\theta_3=-k$ \\ \hline
 IV   & $D=3$, $b_2 > 1$ & $\theta_0 = k, \; \theta_1 = \sqrt{b_2}, \; \theta_2 = -\sqrt{b_2}, \;\theta_3=-k$ \\ \hline
 V   &$D=3$, $b_2 > 1$ & $\theta_0 = k, \; \theta_1 = -\sqrt{b_2}, \; \theta_2 = \sqrt{b_2}, \; \theta_3=-k$ \\ \hline
\end{tabular}
\label{tabela1}
\caption{Special cases}
\end{center}
\end{table}

\begin{lemma}
\label{lem:except}
With reference to Notation \ref{blank1}, assume the $Q$-polynomial structure is listed in Table $1$. Then 
the corresponding dual eigenvalue sequence $\{\ths_i\}_{i=0}^D$ is given 
in the table below.
\renewcommand\arraystretch{2}
\begin{center}
\begin{tabular}[b]{|c|c|}
 \hline 
 {\rm Case}            & {\rm dual eigenvalue sequence}                                                                          \\ \hline \hline
 {\rm I}               & $\ths_i = D-2i \quad (0 \le i \le D)$                                                                   \\ \hline
 {\rm II}              & $\ths_i = 2(D-i)^2 - D \quad (0 \le i \le D)$                                                           \\ \hline
 {\rm III}             & $\ths_i = \theta_i \quad (0 \le i \le 3)$                                                               \\ \hline
 {\rm IV, V}           & $\ths_0 = \frac{k(k-1)}{c_2}, \; \ths_1 = \frac{\theta_1(k-1)}{c_2}, \; \ths_2 = -1, \; \ths_3=-\frac{k}{\theta_1}$     \\ \hline
\end{tabular}
\end{center}
\end{lemma}
\proof
The $\{\ths_i\}_{i=0}^D$ are computed using \eqref{eq:de} with $\theta=\theta_1$, once $\ths_0$ is known.
Recall that $\ths_0$ is the rank of $E_1$. 
In Case I the rank of $E_1$ is $D$ by \cite[Theorem 9.2.1]{BCN}.
In Case II the rank of $E_1$ is $2D^2-D$ by \cite[p. 264]{BCN}.
In Case III  the rank of $E_1$ is $k$ by \cite[p. 432]{BCN}.
In Cases IV and V the rank of $E_1$ is $k(k-1)/c_2$ by \cite[p. 432]{BCN}.
The result follows. \qed

\begin{lemma}
\label{lem:except1}
With reference to Notation \ref{blank1}, assume the $Q$-polynomial structure is listed in Table $1$. Then 
$\beta$, $\gamma^*$, $\varrho$, $\varrho^*$ are given in the table below.
\renewcommand\arraystretch{2}
\begin{center}
\begin{tabular}{|c|c|c|c|c|}
 \hline 
 {\rm Case}      & $\beta$                       & $\gamma^*$                      & $\varrho$              & $\varrho^*$                     \\ \hline \hline
 {\rm I}         & $2$                           & $0$                             & $4$                    & $4$                             \\ \hline
 {\rm II}        & $2$                           & $4$                             & $16$                   & $4(2D-1)$                       \\ \hline
 {\rm III}       & $k-1$                         & $0$                             & $k+1$                  & $k+1$                           \\ \hline
 {\rm IV, V}     & $\frac{k-\theta_1}{\theta_1}$ & $\frac{(k-1)\theta_1-c_2}{c_2}$ & $\theta_1(\theta_1+k)$ & $\frac{(k-1)(k+\theta_1)}{c_2}$ \\ \hline
\end{tabular}
\end{center}
\end{lemma}
\proof
\medskip \noindent
Use Lemma \ref{lem:fundamental1} and Lemma \ref{lem:except}. \qed 

\noindent
We have now completed our description of the special cases.

\begin{lemma}
\label{lem:caug}
{\rm (\cite[Lemma 3.2, Lemma 3.3]{Caug3})}
With reference to Notation \ref{blank1}, assume the $Q$-polynomial structure is not listed in Table $1$ and $\beta \ne -2$.
Then there exist $q,s^* \in \RR$ such that the following {\rm (i)--(iii)} hold.
\begin{itemize}
\item[{\rm (i)}]    $|q| > 1 \;$, $s^* q^i \ne 1 \;(2 \le i \le 2D+1)$,
\item[{\rm (ii)}]   $b_0=h(q^D-1)=c_D$,
                    $$
                      b_i=\frac{h(q^D-q^i)(1-s^*q^{i+1})}{1-s^*q^{2i+1}}, \quad c_i=\frac{h(q^i-1)(1-s^*q^{D+i+1})}{1-s^*q^{2i+1}} \qquad (1 \le i \le D-1),
                    $$
\item[{\rm (iii)}]  $\theta_i=h(q^{D-i}-q^{i}), \;$ $\ths_i=\ths_0 + h^*(1-q^i)(1-s^* q^{i+1})q^{-i} \quad (0 \le i \le D)$, 
\end{itemize}
where
$$
  h=\frac{1-s^* q^3}{(q-1)(1-s^* q^{D+2})}, \; \;
  h^*=\frac{(q^D+q^2)(q^D+q)}{q(q^2-1)(1-s^* q^{2D})}, \; \;
  \ths_0=\frac{h^*(q^D-1)(1-s^* q^2)}{q(q^{D-1}+1)}.
$$
\end{lemma}

\begin{note} 
{\rm With reference to Notation \ref{blank1}, assume the $Q$-polynomial structure is not listed in Table $1$ and $\beta \ne -2$.
Then by \cite[Corollary 6.7]{Caug3}
the scalar $s^*$ from Lemma \ref{lem:caug} is zero provided $D \ge 12$.}
\end{note}

\begin{lemma}
\label{lem:caug_par}
With reference to Notation \ref{blank1}, assume the $Q$-polynomial structure is not listed in Table $1$ and $\beta \ne -2$. 
Then
\begin{align*}
  \beta     & = q+q^{-1}, \\
  \gamma^*  & = \frac{(q-1)(q^{D-2}+1)(1 + s^* q^{D+1})}{1-s^*q^{2D}}, \\
  \varrho   & = \frac{q^{D-2}(q+1)^2(1-s^*q^3)^2}{(1-s^*q^{D+2})^2}, \\
  \varrho^* & = \frac{q(q^{D-2}+1)^2(1-s^*q^2)}{1-s^* q^{2D}},
\end{align*}
where $q, s^*$ are from Lemma \ref{lem:caug}.
\end{lemma}
\proof
Use Lemma \ref{lem:fundamental1} and Lemma \ref{lem:caug}. \qed

\medskip \noindent
In Lemma \ref{lem:caug}(i) we cited some inequalities involving $q$ and $s^*$.
We now prove one more inequality involving $q$ and $s^*$. 

\begin{lemma}
\label{cor:sq}
With reference to Notation \ref{blank1}, assume the $Q$-polynomial structure is not listed in Table $1$ and $\beta \ne -2$. 
Then the scalars $q$ and $s^*$ from Lemma \ref{lem:caug} satisfy $s^*q \ne 1$.
\end{lemma}
\proof
We assume $s^*q=1$ and get a contradiction. Recall $q \in \RR$ and $|q| > 1$.
By \cite[Theorem 15.6(ii)]{Caug} the scalar
$$
  \frac{q(q^{D-1}-1)}{q^{D+1}-1}
$$
is nonnegative. The factors $q^{D-1}-1$ and $q^{D+1}-1$ have the same sign, since $D-1$ and $D+1$ have the same parity.
Therefore $q > 0$. By these comments $q > 1$. 
By \cite[Theorem 15.6(iii)]{Caug} the scalar
\begin{equation}
\label{eq_mod}
  \frac{(q^D-1)(q^D-q)(1-q^3)(1+q^D)}{q(q^2-1)(1-q^{D+1})(1-q^{2D-1})}
\end{equation}
is nonnegative.
Since $q > 1$ the expression \eqref{eq_mod} is negative, for a contradiction. \qed

\section{The main result}
\label{sec:qs1}

Recall our $Q$-polynomial structure from Notation \ref{blank1}.
We are now ready to determine whether the corresponding $R/L$ dependency structure is  uniform or strongly uniform.
We begin with some computations involving the matrix $U$ from Definition \ref{def:ef}.

\renewcommand\arraystretch{2}
\begin{proposition}
\label{thm:par}
With reference to Notation \ref{blank1} and Definition \ref{def:ef}, the scalars 
$\{e_i^-\}_{i=2}^D$, $\{e_i^+\}_{i=1}^{D-1}$, $\{f_i\}_{i=1}^D$
are given in the following table:
\begin{center}
\begin{tabular}{|c|c|c|c|}
 \hline 
 {\rm Case}     & $e_i^-$                                 
 & $e_i^+$        & $f_i$                                     
 \\ \hline \hline
 {\rm I}        & $-\frac{1}{2}$                                       & $-\frac{1}{2}$                                      & $1$                                                \\ \hline
 {\rm II}       & $\frac{i-D-2}{2D-2i+1}$           & $\frac{i-D+1}{2D-2i+1}$                            & $4$                                                \\ \hline
 {\rm III}      & $e_2^-=\frac{2-k}{2}, \; e_3^-=\frac{1}{1-k}$        & $e_1^+=\frac{1}{1-k}, \; e_2^+=\frac{2-k}{2}$       & $1$                                                \\ \hline
 {\rm IV, V}    & $e_2^-=\frac{\theta_1-k+1}{\theta_1+1}, \;
                   e_3^-=-\frac{\theta_1^2}{c_2}$                      & $e_1^+=\frac{1}{1-k}, \; 
                                                                          e_2^+=\frac{\theta_1-c_2}{\theta_1(\theta_1+1)}$   & $\theta_1^2$                                       \\ \hline
 {\rm other}    & $-\frac{q^2(1 - s^* q^{2i-3})}{(q+1)(1-s^* q^{2i})}$ & $-\frac{1-s^* q^{2i+3}}{q(q+1)(1-s^* q^{2i})}$      & $\frac{q^{D-1}(1-s^* q^3)^2}{(1 - s^* q^{D+2})^2}$ \\ \hline 
\end{tabular}
\end{center}
The scalars $q,s^*$ are from Lemma \ref{lem:caug}.
\end{proposition}
\proof
For Cases I--V use Definition \ref{def:ef}, Lemma \ref{lem:except} and Lemma \ref{lem:except1}.
For the remaining case use Definition \ref{def:ef}, Lemma \ref{lem:caug} and Lemma \ref{lem:caug_par}. \qed

\medskip \noindent
With reference to Notation \ref{blank1} assume for the moment that $D=3$.
For an integer $t \ge 2$ the following are equivalent: (i) $\G$ is of McFarland type with parameters $(1,t)$;
(ii) the intersection numbers of $\G$ satisfy $k=t(t+1)$ and $c_2=t$.
Assume that (i), (ii) hold. 
Then $b_1 = t^2+t-1$, $b_2 = t^2$, $c_3=t(t+1)$. 
Moreover the eigenvalue $\theta_1$ is either $t$ or $-t$. 
The case $\theta_1=t$ is contained in Case IV. 
We call this situation Case IV'. Let us examine Case IV' in more detail.

\begin{lemma}
\label{lem:exception}
With reference to Notation \ref{blank1}, assume the $Q$-polynomial structure is in Case IV'.
Then the following {\rm (i)--(iii)} hold.
\begin{itemize}
\item[{\rm (i)}]
The eigenvalues $\{\theta_i\}_{i=0}^3$ are
$$
  \theta_0=t(t+1), \quad \theta_1=t, \quad \theta_2=-t, \quad \theta_3=-t(t+1).
$$
\item[{\rm (ii)}]
The dual eigenvalues $\{\ths_i\}_{i=0}^3$ are
$$
  \ths_0=(t+1)(t^2+t-1), \quad \ths_1=t^2+t-1, \quad \ths_2=-1, \quad \ths_3=-t-1.
$$
\item[{\rm (iii)}]
The parameters $\beta, \gamma^*, \varrho, \varrho^*$ from Lemma \ref{lem:fundamental1} are 
$$
  \beta=t, \quad \gamma^*=t^2+t-2, \quad \varrho=t^2(t+2), \quad \varrho^*=t^3+3t^2+t-2.
$$
\end{itemize}
\end{lemma}
\proof
(i) Immediate from Table 1.

\smallskip \noindent
(ii) Immediate from Lemma \ref{lem:except}.

\smallskip \noindent
(iii) Immediate from Lemma \ref{lem:except1}. \qed

\begin{lemma}
\label{lem:exception1}
With reference to Notation \ref{blank1}, assume the $Q$-polynomial structure is in Case IV'.
Then the following {\rm (i)--(iii)} hold.
\begin{itemize}
\item[{\rm (i)}]
$e_2^- = t-1$ and $e_3^-=t$.
\item[{\rm (ii)}]
$e_1^+ = (t^2+t-1)^{-1}$ and $e_2^+=0$.
\item[{\rm (iii)}]
$f_i=t^2$ for $1 \le i \le 3$.
\end{itemize}
\end{lemma}
\proof
From Case IV of the table of Proposition \ref{thm:par}, using $k=t(t+1)$, $c_2=t$, and $\theta_1=t$. \qed

\begin{corollary}
\label{cor:par}
With reference to Notation \ref{blank1} and Definition \ref{def:ef} the following {\rm (i)--(iii)} hold.
\begin{itemize}
  \item[{\rm (i)}]
  Assume the $Q$-polynomial structure is in Case II. Then $e_{D-1}^+ =0$.
  \item[{\rm (ii)}]
  Assume the $Q$-polynomial structure is in Case IV'. Then $e_2^+=0$.
  \item[{\rm (iii)}]
  For all other cases 
  $e_i^- \ne 0$ for $2 \le i \le D$ and $e_i^+ \ne 0$ for $1 \le i \le D-1$.
\end{itemize}
\end{corollary}
\proof
(i) Immediate from Case II of the table in Proposition \ref{thm:par}.

\smallskip \noindent
(ii) Immediate from Lemma \ref{lem:exception1}(ii).

\smallskip \noindent
(iii) Immediate from Proposition \ref{thm:par},
using Lemma \ref{lem:caug}(i) and Lemma \ref{cor:sq}. \qed

\medskip \noindent
We recall a result from linear algebra.

\begin{lemma}
\label{lem:det}
{\rm (\cite[p. 29]{HJ})}
Pick an integer $d \ge 3$ and let $B = (B_{ij})_{1 \le i,j \le d}$ denote a tridiagonal matrix. 
Then
$$
  \det(B) = B_{dd} \det\big( (B_{ij})_{1 \le i,j \le d-1} \big) - B_{d-1,d} B_{d,d-1} \det\big( (B_{ij})_{1 \le i,j \le d-2} \big).
$$
\end{lemma}
Recall the principal submatrices $(e_{ij})_{r \le i,j \le p}$  from the beginning of Section \ref{sec:uniform}.

\begin{proposition}
\label{thm:det}
With reference to Notation \ref{blank1} and Definition \ref{def:ef}, for $1 \le r \le p \le D$
the determinant of $(e_{ij})_{r \le i,j \le p}$ 
is given in the following table:
\begin{center}
\begin{tabular}{|c|c|}
 \hline 
 {\rm Case}   & {\rm determinant of} 
                  $(e_{ij})_{r \le i,j \le p}$ 
 \\ \hline \hline
{\rm I}      & $\frac{p-r+2}{2^{p-r+1}}$                                                                                     \\ \hline
 {\rm II}     & $\frac{(p-r+2)(2D-r-p+1)(D-p+1)_{p-r}}{2^{p-r+2}(D-p+1/2)_{p-r+1}}$                                           \\ \hline
 {\rm III}    & $1$ if $p=r$, \qquad $\frac{k}{2(k-1)}$ if $p=r+1$, \qquad $\frac{1}{k-1}$ if $p=r+2$                         \\ \hline
 {\rm IV, V}  & \begin{tabular}{l}
                 $1$ if $p=r$ \qquad \qquad \qquad \;
                 $\frac{k \theta_1}{(k-1)(\theta_1+1)}$ if $(r,p)=(1,2)$ \\
                 $\frac{k}{c_2(\theta_1+1)}$ if $(r,p)=(2,3)$ \qquad
                 $\frac{\theta_1(k - \theta_1)}{(k-1)c_2}$ if $(r,p)=(1,3)$
                \end{tabular}                                                                                                 \\ \hline
 {\rm other}  & $\frac{(q^{p-r+2}-1) (1-s^*q^{p+r})(s^* q^{2r+1} ; q^2)_{p-r}}{(q+1)^{p-r+1}(q-1)(s^* q^{2r} ; q^2)_{p-r+1}}$ \\
 \hline 
\end{tabular}
\end{center}
We are using the notation 
$$
  (a)_n = a(a+1) \cdots (a+n-1),
$$
$$
  (a;q)_n = (1-a)(1-qa) \cdots (1-q^{n-1}a).
$$
\end{proposition}
\proof
For Cases III--V the result follows from a straightforward computation using Proposition \ref{thm:par}.
For the other cases use Proposition \ref{thm:par}, Lemma \ref{lem:det} and induction on $p-r$.  \qed

\begin{corollary}
\label{cor:nonsing}
With reference to Notation \ref{blank1} and Definition \ref{def:ef}, 
for $1 \le r \le p \le D$ the principle submatrix 
$(e_{ij})_{r \le i,j \le p}$ 
is nonsingular.
\end{corollary}
\proof
Immediate from Lemma \ref{lem:caug}(i) and Proposition \ref{thm:det}. \qed 

\begin{proposition}
\label{main:pro8}
With reference to Notation \ref{blank1} assume $\beta \ne -2$.
For Case II and Case IV' the corresponding 
$R/L$-dependency structure is
uniform but not strongly uniform.
In all other cases the corresponding 
$R/L$-dependency structure is strongly uniform.
\end{proposition}
\proof
Immediate from Proposition \ref{main:prop5}, 
Corollary \ref{cor:par} and Corollary \ref{cor:nonsing}. \qed

\begin{theorem}
\label{main:thm}
Let $\G$ denote a bipartite distance-regular graph with diameter $D \ge 3$ and valency $k \ge 3$.
Fix a vertex $x$ and let $L$ (resp. $R$) denote the corresponding lowering (resp. raising) matrix
from Definition \ref{def:LR}. Let $\{\theta_i\}_{i=0}^D$ denote a $Q$-polynomial ordering of the
eigenvalues of $\G$.
Consider the following cases:
\begin{itemize}
\item[{\rm (i)}]
$\G$ is the hypercube $H(D,2)$ with $D$ even and $\theta_i=(-1)^i (D-2i)$ for $0 \le i \le D$; 
\item[{\rm (ii)}]
$\G$ is the antipodal quotient $\overline{H}(2D,2)$ and $\theta_i = 2D - 4i$ for $0 \le i \le D$;
\item[{\rm (iii)}]
$D=3$ and $\G$ is of McFarland type with parameters $(1,t)$ for some integer $t \ge 2$, and 
$\theta_0, \theta_1, \theta_2, \theta_3$ are $t(t+1), t, -t, -t(t+1)$ respectively.
\end{itemize}
In Case {\rm (i)} the corresponding 
$R/L$-dependency structure is not uniform.
In Cases {\rm (ii), (iii) } this structure
is uniform but not strongly uniform.
In all other cases this structure is strongly uniform.
\end{theorem}
\proof Immediate from Propositions \ref{pro:beta-2} and \ref{main:pro8}. \qed


\end{document}